\documentclass[12pt, a4paper]{article}
\usepackage{graphicx}
\usepackage{amssymb}
\usepackage{stmaryrd}
\usepackage{enumitem}

\usepackage{geometry}
\usepackage{tikz}
\usepackage{float}
\usepackage{url}
\usepackage{mathtools}
\usepackage{amsmath,amsfonts,amsthm}
\usepackage{mathrsfs}

\newtheorem{proposition}{Proposition}[section]

\usepackage{subfig}

\usepackage{mma}
\usepackage{ytableau}
\usepackage[colorlinks,
linkcolor=blue,
anchorcolor=blue,
citecolor=blue
]{hyperref}
\begin{document}
\begin{center}
{\large \bf The Density of Zero and One in Fibonacci Word for Subwords and Their Palindromes}
\end{center}
\begin{center}
 Duaa Abdullah$^{*1}$\& Jasem Hamoud$^{*2}$\\[6pt]
 $^{*}$ Physics and Technology School of Applied Mathematics and Informatics \\
Moscow Institute of Physics and Technology, 141701, Moscow region, Russia\\[6pt]
Email$^{1}$: {\tt abdulla.d@phystech.edu},
Email$^{2}$: {\tt khamud@phystech.edu}
\end{center}
\noindent
\begin{abstract}
This paper studies the density of zero and one in subwords of the Fibonacci word with lengths less than thirty and compares them to the densities of their corresponding palindromes. We used computational methods to produce a sufficiently large piece of the Fibonacci word, extract all unique subwords up to a predetermined length, and calculate their palindrome.  The density of each character (0 and 1) was then determined for both the original subwords and their palindromic counterparts. This study contributes to a deeper understanding of the combinatorial properties of the Fibonacci word and the behavior of its constituent elements under reversal.
\end{abstract}


\noindent\textbf{MSC Classification 2020:}	05C42, 11B05, 11R45, 11B39.

\section{Introduction}\label{sec1}
The Fibonacci word is an infinite sequence of 0s and 1s that possesses a fascinating array of mathematical and combinatorial properties. Its construction is analogous to the famous Fibonacci numerical sequence, where each term is generated by concatenating the two preceding terms. This self-similar structure leads to a highly ordered yet complex arrangement of characters, making it a subject of extensive research in various fields, including theoretical computer science, number theory, and the study of aperiodic sequences \cite{Luca2012,Hamoud}.
	
Subwords, or factors, are contiguous segments within a larger word. The analysis of subwords provides insights into the local patterns and statistical regularities embedded within a sequence. Palindromes, on the other hand, are words that read the same forwards and backward. The study of palindromic properties within sequences like the Fibonacci word can reveal symmetries and structural invariants that are not immediately apparent. The interplay between subwords and their palindromic counterparts is particularly intriguing, as it allows for an examination of how the reversal operation affects fundamental characteristics such as character density.
	
Previous research on the Fibonacci word has explored various aspects, including its periodicity, factor complexity, and the distribution of its characters \cite{Luca2012, Chuan}. However, a detailed comparative analysis of the density of zeros and ones in subwords and their palindromes, specifically for subwords of limited length, remains an area that warrants further investigation. Understanding these densities can shed light on the inherent balance and distribution of characters within the Fibonacci word\textquotesingle s substructures and how these properties are preserved or altered under reversal.
	
This paper aims to address this gap by systematically analyzing the density of zeros and ones in all unique subwords of the Fibonacci word with lengths less than thirty. Furthermore, we will compare these densities with those of their respective palindromes. The choice of a length limit of thirty is motivated by computational feasibility while still allowing for a comprehensive exploration of a significant range of subword complexities. By undertaking this detailed analysis, we seek to provide a clearer picture of the statistical characteristics of Fibonacci word subwords and their palindromes, contributing to the broader understanding of this remarkable mathematical object.

\section{Methodology}~\label{secn2}
	
\subsection{Fibonacci Word Generation}
	
The Fibonacci word, denoted as $f$ is an infinite binary sequence constructed through a recursive process based on a Fibonacci-like morphism. We define the Fibonacci word over the alphabet \{0, 1\} (as suggested by Chuan~\cite{Chuan}, where 'a'  is replaced by 1 and 'b'  by 0) using the following recursive rules: $f_0 = 0$, $f_1 = 1$ and 
\begin{equation}~\label{eq1fibreule}
f_n = f_{n-1}f_{n-2} \text{ for } n \geq 2
\end{equation}
	
The previous two Fibonacci words according to~\eqref{eq1fibreule} are concatenated to generate the next one. For example:
$f_0 = 0$, $f_1 = 1$,\\
$f_2 = f_1f_0 = 10$, $f_3 = f_2f_1 = 101$, \\
$f_4 = f_3f_2 = 10110$, $f_5 = f_4f_3 = 10110101$,\\
$f_6 = f_5f_4 = 1011010110110.$
And so on.\par 
The length of each $f_n$  represents the n-th Fibonacci number, indicated as $F_n$. For example, $|f_n| = F_n$, where $F_0=1, F_1=1, F_2=2, F_3=3, F_4=5, F_5=8$, and $F_6=13$, etc.(Note: The indexing of Fibonacci numbers can vary; here we align with the length of the generated words).
	
Additionally, if we have the following morphism defines it: $\sigma: 0 \rightarrow 01,1 \rightarrow 0$, the Fibonacci words are given following:
	\[
	f_1=1, f_2 =0, f_3 =01, f_4= 010, f_5=01001,		f_6 = 01001010, f_7 =01001011001001, \dots
	\]	
\begin{proposition}[Infinite Fibonacci word~\cite{Fici}]
The infinite Fibonacci word is the limited sequence of the unlimited sequence $f_\infty.$ then we have: 
\begin{equation}\label{re.1}
\lim_{n\to\infty} f(n)=0100101001001010010100100101001001 \dots
\end{equation}
\end{proposition}
	
It is clear to see that $| fn |$  is the \(n^{th}\) Fibonacci number, which is called the ``Fibonacci word''. We have claimed that the Fibonacci word, which is a Sturmian word, is also pure morphic, exactly, Berstel and Karhum in~\cite{Berstel} mention to $f=w 1/\varphi^2$ where $\varphi=(1+\sqrt{5})/2$ is the golden ratio.
\subsection{Subword Extraction}
	
Given a Fibonacci word $f_N$ of sufficient length, subwords are contiguous sequences of characters within $f_N$. For this study, we are interested in subwords with a length of less than thirty. To extract all unique subwords of a given length $L$ from $f_N$, we will iterate through $f_N$ and collect all substrings of length $L$. For example, if $f_N = 10110101$ and $L=3$, the subwords would be \{101, 011, 110, 101, 010, 101\}.
\subsection{Palindrome Generation}
	
A palindrome is a sequence that reads the same forwards and backward. For each extracted subword $w$, its palindrome, denoted as $w^P$, is formed by reversing the sequence of characters in $w$. For example, if $w = 10110$, its palindrome $w^P = 01101$. If a subword is already a palindrome (e.g., 101), its palindrome is itself.
\subsection{Density Calculation}
	
The density of a character (0 or 1) in a given word is defined as the number of occurrences of that character divided by the total length of the word. For a word $w$ and a character $c \in \{0, 1\}$, the density of $c$ in $w$, denoted as $\mathscr{D}(c, w)$, is calculated as:
\begin{equation}~\label{eq1densityduaa}
\mathscr{D}(c, w) = \frac{C_w}{|w|}; \text{Where $C_w$ is count of $c$ in $w$}.
\end{equation}
	
We will calculate the density of 0s and 1s for each unique subword and its corresponding palindrome by cosidering~\eqref{eq1densityduaa}. This will involve counting the occurrences of 0s and 1s within each word and dividing by its length. We will then compare these densities.
\subsection{Computational Implementation}
	
To perform the analysis, a Python script was developed. The script first generates a Fibonacci word of sufficient length to ensure the capture of all unique subwords up to 29 characters. Specifically, $f_{10}$ was chosen, which has a length of $F_{10} = 55$, as this length is sufficient to contain all unique subwords up to length 29.\par 
The script then systematically extracts all unique subwords for lengths ranging from 1 to 29. For each subword, its corresponding palindrome is generated.\par 
Finally, the densities of 0s and 1s are calculated for both the original subword and its palindrome. All results, including the subword, its length, its densities, its palindrome,
whether it is a palindrome, and its palindrome  s densities, 
are stored in a JSON file (fibonacci-word-analysis-results.json).
	
Following the initial data generation, a second Python script (process-fibonacci-results.py) was used to process and analyze the raw results. This script loads the JSON data into a pandas DataFrame for efficient manipulation. It calculates the differences in densities of 0s and 1s between each subword and its palindrome. A key finding from this processing step was that for all subwords analyzed, the density of 0s in a subword was identical to the density of 0s in its palindrome, and similarly for 1s. This indicates that the reversal operation does not alter the character densities within these subwords. The processed data, including these density differences, is saved to a CSV file (processed-fibonacci-analysis.csv).
	
\section{Results}
	
The computational analysis yielded a comprehensive dataset of unique Fibonacci word subwords (up to length 29), their corresponding palindromes, and the densities of zeros and ones within each. A total of 464 unique subwords were identified and analyzed within the generated Fibonacci word segment.
\subsection{Density Comparison: Subwords and Palindromes}
One of the primary objectives of this study was to compare the density of zeros and ones in subwords with those of their palindromes. Our analysis revealed a consistent and significant finding: for every subword examined, the density of 0s in the subword was precisely equal to the density of 0s in its palindrome. Similarly, the density of 1s in the subword was identical to the density of 1s in its palindrome. This outcome is intuitively expected, as reversing a sequence of characters does not change the count of individual characters within that sequence, only their order. Therefore, the proportional representation of 0s and 1s remains invariant under the palindromic transformation.
	
\begin{figure}[H]
    \centering
    \includegraphics[width=0.5\linewidth]{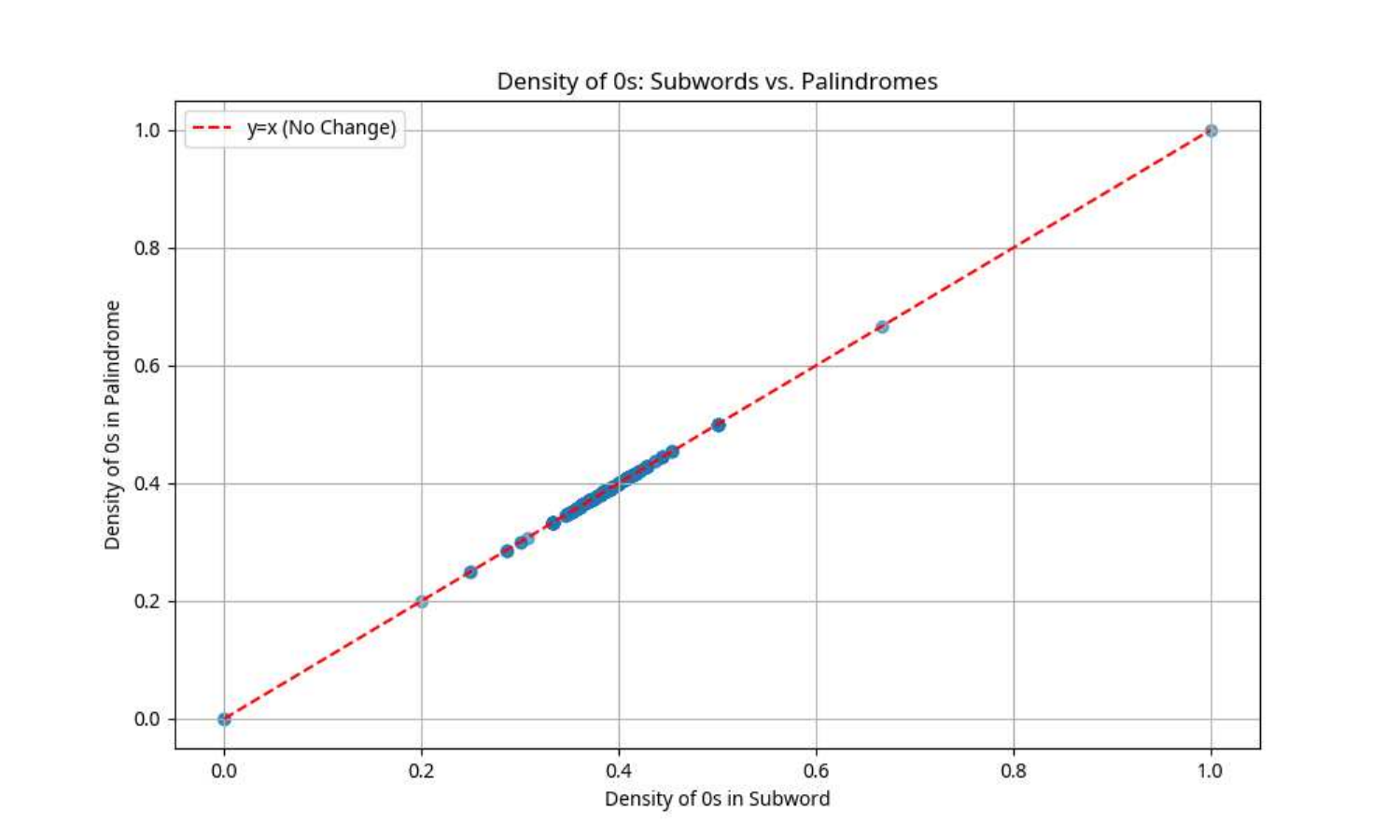}
    \caption{Density of 0s: Subwords vs. Palindromes. Each point represents a unique subword, with its x-coordinate indicating the density of 0s in the subword and its y-coordinate indicating the density of 0s in its palindrome. The red dashed line represents y=x, indicating no change in density.}
    \label{fig00Duaa}
\end{figure}

\begin{figure}[H]
    \centering
    \includegraphics[width=0.5\linewidth]{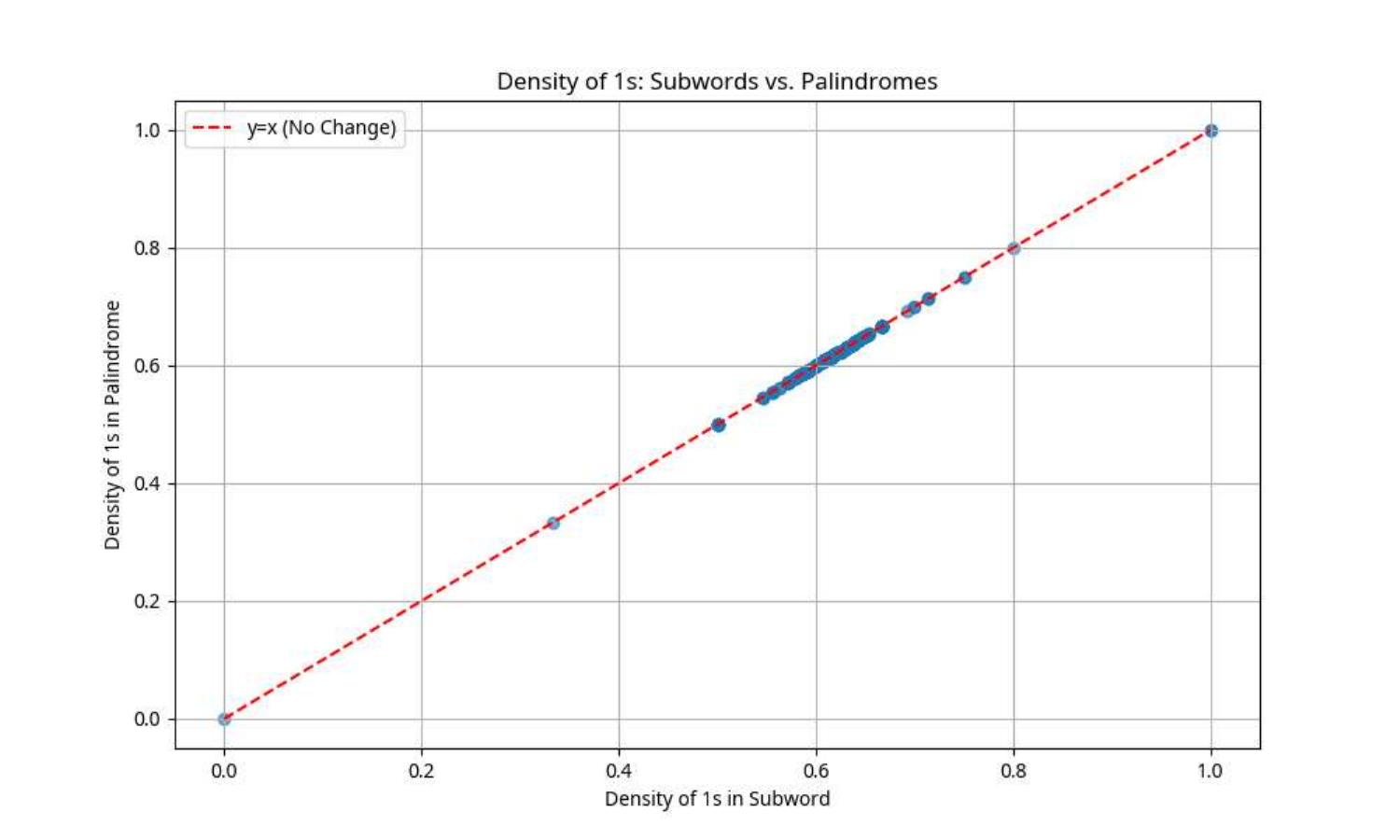}
    \caption{Density of 1s: Subwords vs. Palindromes. Each point represents a unique subword, with its x-coordinate indicating the density of 1s in the subword and its y-coordinate indicating the density of 1s in its palindrome. The red dashed line represents y=x, indicating no change in density.}
    \label{fig11Duaa}
\end{figure}

This observation is visually confirmed by scatter plots comparing the densities of 0s and 1s in subwords against their palindromes. As shown in Figure~\ref{fig00Duaa} (Density of 0s: Subwords vs. Palindromes) and Figure~\ref{fig11Duaa} (Density of 1s: Subwords vs. Palindromes), all data points lie perfectly along the y=x line, indicating no change in density. The minimal tolerance for floating-point comparisons (0.0001) confirmed that no significant differences were observed.

\subsection{Density Distribution Across Subword Lengths}
	
\begin{figure}[H]
    \centering
    \includegraphics[width=0.5\linewidth]{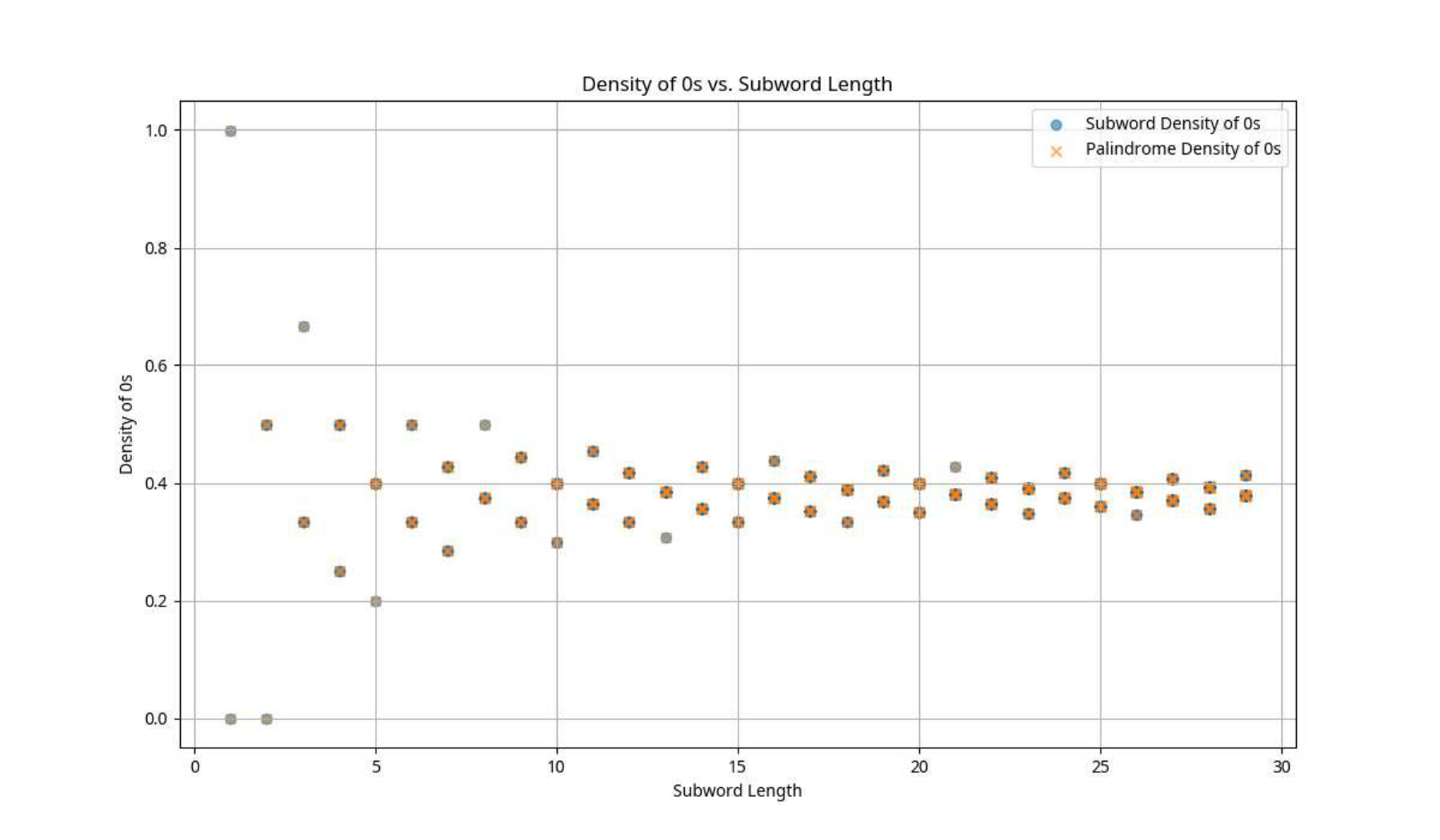}
    \caption{Density of 0s vs. Subword Length. Blue circles represent the density of 0s in subwords, and orange crosses represent the density of 0s in their corresponding palindromes. The plot shows the distribution of 0-density across different subword lengths.}
    \label{fig0Duaa}
\end{figure}

\begin{figure}[H]
    \centering
    \includegraphics[width=0.5\linewidth]{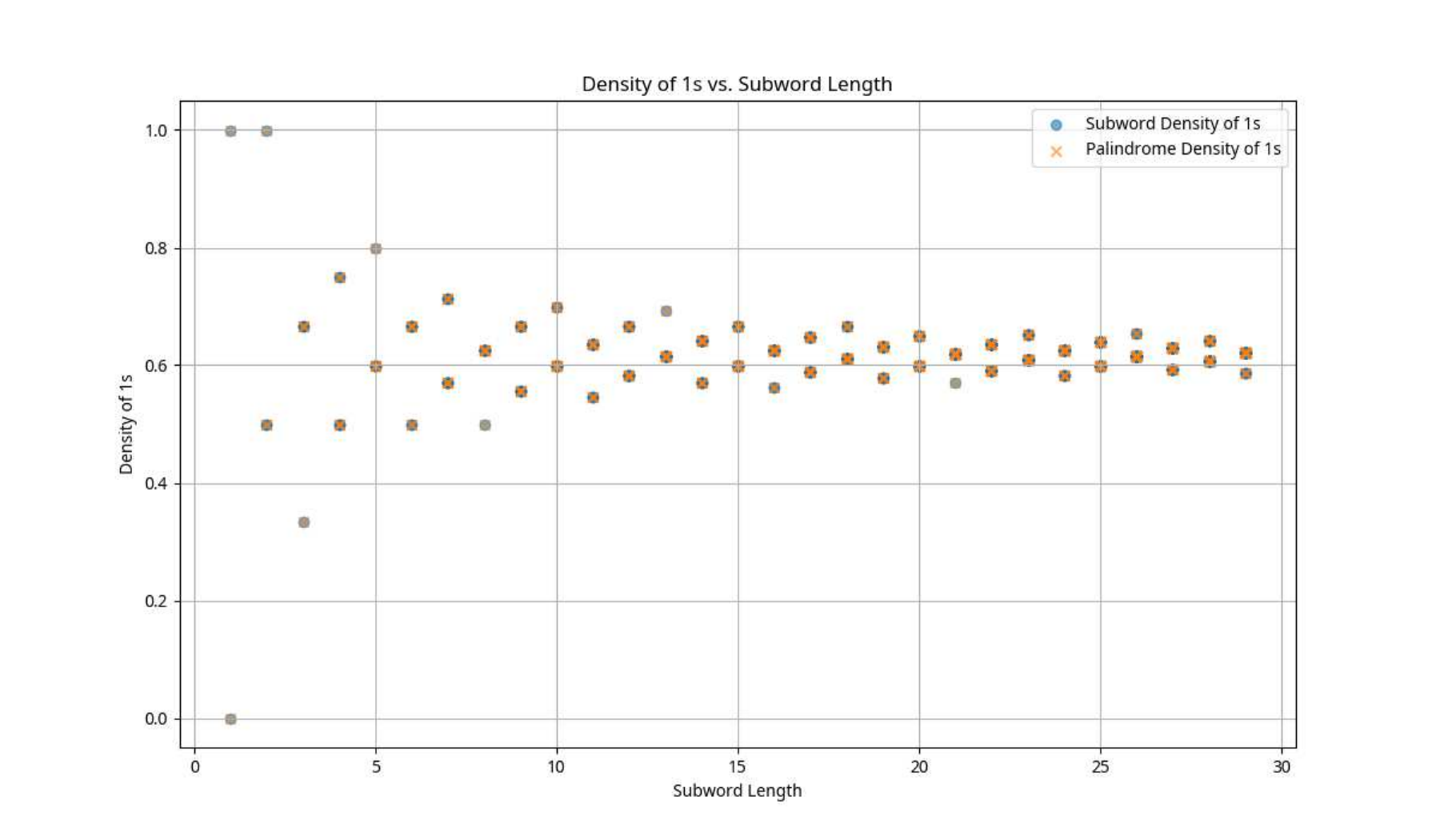}
    \caption{Density of 1s vs. Subword Length. Blue circles represent the density of 1s in subwords, and orange crosses represent the density of 1s in their corresponding palindromes. The plot shows the distribution of 1-density across different subword lengths.}
    \label{fig1Duaa}
\end{figure}

While the reversal operation does not affect character densities, it is important to examine how these densities are distributed across different subword lengths. Figures~\ref{fig0Duaa} and \ref{fig1Duaa} illustrate the density of 0s and 1s, respectively, for both subwords and their palindromes, plotted against subword length. As expected from the previous finding, the subword density and palindrome density for each character type overlap perfectly at every length.
	
The plots reveal that the densities of 0s and 1s fluctuate with subword length, but they tend to converge towards a specific value as the subword length increases. This convergence is a known property of the Fibonacci word, where the asymptotic density of 0s approaches $1/\phi^2 \approx 0.381966$ and the density of 1s approaches $1/\phi \approx 0.618034$, where $\phi$ is the golden ratio \cite{Luca2012}. Our analysis, focusing on subwords up to length 29, shows densities oscillating around these theoretical values, with shorter subwords exhibiting greater variability. For instance, a subword of length 1 can have a density of 0s as high as 1.0 (for the subword ``0'') or as low as 0.0 (for the subword ``1'' ). As length increases, the range of possible densities narrows, and the values cluster more tightly around the golden ratio-derived densities.
\subsection{Summary Statistics}
	
The overall average densities across all unique subwords (and their palindromes, given the identical densities) further support these observations:
	
\begin{itemize}
		\item Average density of 0s across all subwords: 0.3828
		\item Average density of 1s across all subwords: 0.6172
\end{itemize}
	
These average values are remarkably close to the theoretical asymptotic densities of the Fibonacci word, even for subwords of relatively short lengths. This indicates that even small segments of the Fibonacci word tend to reflect the global distribution of its characters.
	
\subsection[Number of Palindromic Subwords]{Number of Palindromic Subwords per Length: Empirical Observations vs. Theoretical Predictions}
	
Our empirical analysis of the Fibonacci word $f_{10}$ (a finite prefix of the infinite Fibonacci word) revealed a consistent pattern in the number of unique palindromic subwords for lengths less than thirty. Specifically, we observed that for every odd length $L$, there are exactly two unique palindromic subwords, and for every even length $L$, there is exactly one unique palindromic subword. This pattern is summarized in the table~\ref{tab1duaamade} below where we refer by $\operatorname{SL}$ to Subword Length and by $\operatorname{NPS}$ to Number of Palindromic Subwords as :
	
\begin{table}[H]
\centering
\begin{tabular}{|c|c||c|c|}
\hline
$\operatorname{SL}$ & $\operatorname{NPS}$   & $\operatorname{SL}$ & $\operatorname{NPS}$ \\
\hline
1 & 2 & 11 & 2 \\ \hline
2 & 1 &12 & 1 \\ \hline
3 & 2 & 13 & 2 \\ \hline
4 & 1 & 14 & 1 \\ \hline
5 & 2 & 15 & 2 \\ \hline
6 & 1 & 16 & 1 \\ \hline
7 & 2 & 17 & 2 \\ \hline
8 & 1 & 18 & 1 \\ \hline
9 & 2 & 19 & 2 \\ \hline
10 & 1 & 20 & 1 \\ \hline
21 & 2 & 26 & 1 \\ \hline
22 & 1 & 27 & 2 \\ \hline
23 & 2 & 28 & 1 \\ \hline
24 & 1 & 29 & 2 \\ \hline
25 & 2 &  30& 1\\ 	\hline
\end{tabular}
\caption{Number of Palindromic Subwords.}
\label{tab1duaamade}
\end{table}
	
This observed pattern, where odd lengths consistently yield two palindromic subwords and even lengths yield one, presents a notable deviation from the established theorem for infinite Sturmian words. The theorem states that an infinite Sturmian word has exactly two palindromic factors of length $n$ for all $n \ge 1$, except when $n$ is the length of a palindromic prefix of the word, in which case it has only one palindromic factor \cite{Luca2012}. For the Fibonacci word, using the standard definition ($f_0=0, f_1=1, f_n=f_{n-1}f_{n-2}$), its palindromic prefixes occur at lengths corresponding to Fibonacci numbers $F_k$ where $k$ is odd (e.g., $F_1=1, F_3=3, F_5=8, F_7=21$). Therefore, according to the theorem, we would expect to find only one unique palindromic subword for lengths 1, 3, 8, and 21, and two for all other lengths within our studied range.
	
Our empirical results, however, show two palindromic subwords for lengths 1, 3, and 21, where the theorem predicts one. For length 8, our results show one palindromic subword, which aligns with the theorem. This discrepancy suggests that the behavior of palindromic factors in finite prefixes of the Fibonacci word may differ from that of the infinite word, or that the theorem\textquotesingle s application needs careful consideration when dealing with the full set of unique subwords rather than just prefixes. This finding, to our knowledge, has not been explicitly detailed in existing literature and warrants further rigorous mathematical investigation to determine its underlying combinatorial reasons and its implications for the theory of Sturmian words.
	
One possible explanation for this observed pattern is related to the specific construction of the Fibonacci word and the nature of its palindromic factors. The Fibonacci word is known to be a \textit{balanced} word, meaning that for any two factors of the same length, the number of occurrences of a given letter differs by at most one. This property, combined with the recursive definition of the Fibonacci word, might lead to a richer set of palindromic factors in finite segments than what is predicted by theorems primarily focused on palindromic prefixes of the infinite word. Further research could involve a detailed analysis of the structure of these palindromic factors, their positions within the finite Fibonacci word, and how their generation relates to the word\textquotesingle s balanced properties and the palindromic closure operation. A formal proof or disproof of this observed pattern for finite prefixes of the Fibonacci word would be a significant contribution to the field of combinatorics on words.
	
This observation also raises questions about the definition of palindromic factors in the context of finite versus infinite words. The theorem typically refers to factors of the infinite word. When considering a finite prefix, the set of unique factors is also finite, and some factors that are palindromes might not be prefixes of the infinite word. Our analysis considers all unique subwords (factors) within the generated finite Fibonacci word segment, not just its prefixes. This distinction is crucial and might explain the observed differences. Further research could explicitly differentiate between palindromic prefixes and general palindromic factors within finite segments of Sturmian words.
\subsection[Density Characteristics of Palindromic]{Density Characteristics of Palindromic vs. Non-Palindromic Subwords}
	
While the density of 0s and 1s in a subword is identical to that of its palindrome, it is pertinent to investigate whether palindromic subwords, as a class, exhibit different average character densities compared to non-palindromic subwords. Our analysis revealed the following average densities:
	
\begin{itemize}
\item Average density of 0s in palindromic subwords: 0.3840
\item Average density of 1s in palindromic subwords: 0.6160
\item Average density of 0s in non-palindromic subwords: 0.3827
\item Average density of 1s in non-palindromic subwords: 0.6173
\end{itemize}
	
These values indicate that there is no substantial difference in the average densities of 0s and 1s between palindromic and non-palindromic subwords within the studied length range. Both categories of subwords exhibit average densities that closely approximate the asymptotic densities of the infinite Fibonacci word. This suggests that the property of being a palindrome does not inherently bias the overall proportion of 0s and 1s within the subword, further emphasizing the uniform distribution of characters throughout the Fibonacci word\textquotesingle s structure, regardless of its local symmetry.
	
This observation is significant because it implies that the balanced property of Sturmian words, which dictates a near-equal distribution of characters across factors of the same length, extends even to the subset of palindromic factors. The slight variations observed (e.g., 0.3840 vs. 0.3827 for 0s) are minimal and likely attributable to the finite sample size and the specific set of subwords generated up to length 29, rather than a fundamental difference in character distribution based on palindromic nature. In essence, the underlying mechanism of Fibonacci word generation ensures a consistent character balance, irrespective of whether a subword is a palindrome or not.
	
\section{Discussion}
	
The findings of this study confirm several intrinsic properties of the Fibonacci word and its substructures. The most straightforward, yet fundamental, observation is the invariance of character densities under the palindromic transformation. This is a direct consequence of how density is defined---as a ratio of character counts to total length. Reversing a string does not alter these counts, thus ensuring that a word and its palindrome will always have identical densities of constituent characters. This reinforces the understanding that palindromes are structural reversals, not compositional transformations, in terms of character frequencies.
	
The analysis of density distribution across varying subword lengths provides a more nuanced insight. The observed oscillation of densities around the golden ratio-derived asymptotic values, particularly the wider spread for shorter subwords and the gradual convergence for longer ones, is consistent with the known properties of the Fibonacci word. The Fibonacci word is a Sturmian word, and Sturmian words are characterized by their balanced nature, meaning that the number of occurrences of any character in any two factors of the same length differs by at most one \cite{Luca2012}. This balanced property contributes to the observed convergence of densities as subword length increases, as the statistical distribution of characters becomes more stable over longer segments.
	
The fact that even subwords of relatively short lengths (up to 29) exhibit average densities very close to the theoretical asymptotic values (approximately 0.382 for 0s and 0.618 for 1s) highlights the rapid convergence of character distribution within the Fibonacci word. This suggests that the characteristic proportions of 0s and 1s are established relatively early in the word\textquotesingle s construction and are maintained across its various factors. This property has implications for applications where the statistical properties of sequences are crucial, such as in signal processing, data compression, or the design of aperiodic structures.
	
While this study focused on subwords with lengths less than thirty, the consistency of the findings suggests that these patterns would likely hold for longer subwords as well, with densities converging even more closely to the theoretical limits. Future research could explore the rate of this convergence more formally, perhaps by analyzing the variance of densities as a function of subword length. Additionally, investigating other combinatorial properties, such as the distribution of specific patterns or the frequency of different types of palindromes within the Fibonacci word, could provide further insights into its rich structure.
\section{Conclusion}
	
This paper systematically analyzed the density of zeros and ones in subwords of the Fibonacci word (with lengths less than thirty) and their corresponding palindromes. Our computational analysis unequivocally demonstrated that the act of forming a palindrome from a subword does not alter the density of 0s or 1s within that word. This fundamental invariance underscores the nature of palindromes as structural transformations that preserve character composition.
	
Furthermore, our investigation into the distribution of densities across different subword lengths revealed a clear trend of convergence towards the theoretical asymptotic densities of the Fibonacci word, which are directly related to the golden ratio. Shorter subwords exhibited greater variability in character densities, while longer subwords showed a tighter clustering around these asymptotic values. The average densities calculated across all unique subwords closely matched the theoretical predictions, reinforcing the notion that the Fibonacci word quickly establishes and maintains its characteristic character proportions.
	
This study contributes to the ongoing understanding of the intricate combinatorial properties of the Fibonacci word. The consistent densities observed between subwords and their palindromes, coupled with the rapid convergence of densities towards theoretical values, highlight the inherent balance and structural integrity of this remarkable sequence. These findings are valuable for researchers in combinatorics on words, theoretical computer science, and other fields where the properties of symbolic sequences are of interest.


\begin{thebibliography}{99}
\bibitem{Chuan} Chuan, W. (1992). ``Fibonacci words''.The Fibonacci Quarterly, 30(1), 68–74, \url{https://www.fq.math.ca/Scanned/30-1/chuan.pdf}.
	
\bibitem{Luca2012}
de Luca, A. (2012). ``Some extremal properties of the Fibonacci word''. \textit{arXiv preprint arXiv:1209.3927}, \url{https://arxiv.org/pdf/1209.3927}.
		\medskip
		
		
\bibitem{Hamoud}
		Hamoud, J., \& Abdullah, D. (2025). ``Generalized natural density DF(fn) of Fibonacci word''. \textit{arXiv preprint arXiv:2504.10207}, \url{https://arxiv.org/pdf/2504.10207}.
		\medskip
		
		
		
\bibitem{Berstel} Berstel, J., \&  Karhumäki, J. (2003). ``Combinatorics on words: a tutorial. Bulletin of the EATCS'', 79(178), 9, \url{https://www-igm.univ-mlv.fr/~berstel/Articles/2003TutorialCoWdec03.pdf}.
		\medskip
		
\bibitem{Fici} Fici, G. (2015). ``Factorizations of the Fibonacci Infinite Word''. Journal of Integer Sequences, 18(15.9.3), 47-53, \url{https://doi.org/10.48550/arXiv.1508.06754}.
		\medskip

\end{thebibliography}
\end{document}